\documentclass[english]{amsart}
\usepackage[T1]{fontenc}
\usepackage[latin1]{inputenc}
\usepackage{amssymb}

\makeatletter
 \theoremstyle{plain}
\newtheorem{thm}{Theorem}[section]
  \theoremstyle{plain}
  \newtheorem{lem}[thm]{Lemma}
  \theoremstyle{remark}
  \newtheorem{rem}[thm]{Remark}
  \theoremstyle{remark}
  \newtheorem*{acknowledgement*}{Acknowledgement}

\usepackage{babel}
\makeatother
\begin{document}

\title{On A Limiting Relation Between Ramanujan's Entire Function $A_{q}(z)$
And $\theta$-Function}

\author{Ruiming Zhang}

\date{June 16, 2006}

\curraddr{School of Mathematics\\
Guangxi Normal University\\
Guilin City, Guangxi 541004\\
P. R. China.}

\email{ruimingzhang@yahoo.com}

\keywords{Asymptotics, $\theta$-function, $q$-binomial theorem, $q$-Airy
function, Ramanujan's entire function, discrete Laplace method.}

\subjclass{Primary 33D45. Secondary 33E05.}

\begin{abstract}
We will use a discrete analogue of the classical Laplace method to
show that the main term of the asymptotic expansions of certain entire
functions, including Ramanujan's entire function $A_{q}(z)$, can
be expressed in terms of $\theta$-functions. 
\end{abstract}
\maketitle

\section{Introduction}

Throughout the paper, we assume that \begin{equation}
0<q<1.\label{eq:1.1}\end{equation}
 For any complex number $a$, we define \cite{Andrews4,Gasper,Ismail2}
\begin{equation}
(a;q)_{\infty}=\prod_{k=0}^{\infty}(1-aq^{k})\label{eq:1.2}\end{equation}
 and the $q$-shifted factorial as\begin{equation}
(a;q)_{n}=\frac{(a;q)_{\infty}}{(aq^{n};q)_{\infty}}\label{eq:1.3}\end{equation}
 for any integer $n$. Assume that $|z|<1$, the $q$-Binomial theorem
is \cite{Andrews4,Gasper,Ismail2} \begin{equation}
\frac{(az;q)_{\infty}}{(z;q)_{\infty}}=\sum_{k=0}^{\infty}\frac{(a;q)_{k}}{(q;q)_{k}}z^{k},\label{eq:1.4}\end{equation}
 which defines an analytic function in the region $|z|<1$. Its limiting
case \begin{equation}
(z;q)_{\infty}=\sum_{k=0}^{\infty}\frac{q^{k(k-1)/2}}{(q;q)_{k}}(-z)^{k}\label{eq:1.5}\end{equation}
is one of many $q$-exponential identities. Ramanujan's entire function
$A_{q}(z)$ is defined as \cite{Ismail2}\begin{equation}
A_{q}(z)=\sum_{k=0}^{\infty}\frac{q^{k^{2}}}{(q;q)_{k}}(-z)^{k}.\label{eq:1.6}\end{equation}
 It is known that $A_{q}(z)$ has infinitely many positive zeros and
satisfies the following three term recurrence \begin{equation}
A_{q}(z)-A_{q}(qz)+qzA_{q}(q^{2}z)=0.\label{eq:1.7}\end{equation}
Ramanujan function $A_{q}(z)$, which is also called $q$-Airy function
in the literature, appears repeatedly in Ramanujan's work starting
from the Rogers-Ramanujan identities, where $A_{q}(-1)$ and $A_{q}(-q)$
are expressed as infinite products, \cite{Andrews1}, to properties
of and conjectures about its zeros, \cite{Andrews3,Andrews4,Hayman,Ismail5}.
It is called $q$-Airy function because it appears repeatedly in the
Plancherel-Rotach type asymptotics \cite{Ismail1,Ismail6,Ismail7}
of $q$-orthogonal polynomials, just like classical Airy function
in the classical Plancherel-Rotach asymptotics of classical orthogonal
polynomials \cite{Szego,Ismail2}. However, Ramanujan's $A_{q}(z)$
is not the $q$-analogue of classical Airy functions. Since\begin{equation}
\frac{1-q^{k}}{1-q}\ge kq^{k-1},\label{eq:1.8}\end{equation}
for $k=1,2,...$, then \begin{equation}
\left|\frac{(1-q)^{k}}{(q;q)_{k}}q^{k^{2}}(-z)^{k}\right|\le\frac{(q\left|z\right|)^{k}}{k!}\label{eq:1.9}\end{equation}
 for $k=0,1,\dotsc$, for any complex number $z$, applying Lebesgue's
dominated convergent theorem we have\begin{equation}
\lim_{q\to1}A_{q}((1-q)z)=e^{-z}.\label{eq:1.10}\end{equation}
 We also have obtained the inequality\begin{equation}
\left|A_{q}((1-q)z)\right|\le e^{q|z|}\label{eq:1.11}\end{equation}
for any complex number $z$. For any nonzero complex number $z$,
we define the theta function as \begin{equation}
\theta(z;q)=\sum_{k=-\infty}^{\infty}q^{k^{2}/2}z^{k}.\label{eq:1.12}\end{equation}
 Jacobi's triple product formula says that \cite{Andrews4,Gasper,Ismail2}
\begin{equation}
\sum_{k=-\infty}^{\infty}q^{k^{2}/2}z^{k}=(q,-q^{1/2}z,-q^{1/2}/z;q)_{\infty}.\label{eq:1.13}\end{equation}
 For any real number $x$, then, \begin{equation}
x=\left\lfloor x\right\rfloor +\left\{ x\right\} ,\label{eq:1.14}\end{equation}
where the fractional part of $x$ is $\left\{ x\right\} \in[0,1)$
and $\left\lfloor x\right\rfloor \in\mathbb{Z}$ is the greatest integer
less or equal $x$. The arithmetic function \begin{equation}
\chi(n)=2\left\{ \frac{n}{2}\right\} =n-2\left\lfloor \frac{n}{2}\right\rfloor ,\label{eq:1.15}\end{equation}
 which is the principal character modulo $2$,\begin{equation}
\chi(n)=\begin{cases}
1 & 2\nmid n\\
0 & 2\mid n\end{cases}.\label{eq:1.16}\end{equation}
For any positive real number $t$, we consider the following set \begin{equation}
\mathbb{S}(t)=\left\{ \left\{ nt\right\} :n\in\mathbb{N}\right\} .\label{eq:1.17}\end{equation}
 It is clear that $\mathbb{S}(t)\subset[0,1)$ and it is a finite
set when $t$ is a positive rational number. In this case, for any
$\lambda\in\mathbb{S}(t)$, there are infinitely many positive integers
$n$ and $m$ such that \begin{equation}
nt=m+\lambda,\label{eq:1.18}\end{equation}
where\begin{equation}
m=\left\lfloor nt\right\rfloor .\label{eq:1.19}\end{equation}
If $t$ is a positive irrational number, then $\mathbb{S}(t)$ is
a subset of $(0,1)$ with infinite elements, and it is well-known
that $\mathbb{S}(t)$ is uniformly distributed in $(0,1)$. A theorem
of Chebyshev \cite{Hua} says that given any $\beta\in[0,1)$, there
are infinitely many positive integers $n$ and $m$ such that\begin{equation}
nt=m+\beta+\gamma_{n}\label{eq:1.20}\end{equation}
 with\begin{equation}
|\gamma_{n}|\le\frac{3}{n}.\label{eq:1.21}\end{equation}
 For $n$ large enough, this implies\begin{equation}
m=\left\lfloor nt\right\rfloor .\label{eq:1.22}\end{equation}
 We will also make use of the trivial inequalities \begin{equation}
|e^{x}-1|\le|x|e^{|x|}\label{eq:1.23}\end{equation}
 for any $x\in\mathbb{C}$, and \begin{equation}
e^{-x}\ge1-x\label{eq:1.24}\end{equation}
for $0<x<1$. The following lemma is from \cite{Ismail7}.

\begin{lem}
\label{lem:1}Given any $n\in\mathbb{N}$, if $a>0$, \begin{equation}
\frac{(a;q)_{\infty}}{(a;q)_{n}}=(aq^{n};q)_{\infty}=1+R_{1}(a;n)\label{eq:1.25}\end{equation}
 with\begin{equation}
\left|R_{1}(a;n)\right|\le\frac{(-aq^{2};q)_{\infty}aq^{n}}{1-q}.\label{eq:1.26}\end{equation}
 While for $0<aq<1$, \begin{equation}
\frac{(a;q)_{n}}{(a;q)_{\infty}}=\frac{1}{(aq^{n};q)_{\infty}}=1+R_{2}(a;n)\label{eq:1.27}\end{equation}
 with \begin{equation}
\left|R_{2}(a;n)\right|\le\frac{aq^{n}}{(1-q)(aq;q)_{\infty}}.\label{eq:1.28}\end{equation}
 
\end{lem}
\begin{proof}
It is clear from \eqref{eq:1.5} that\[
R_{1}(a;n)=\sum_{k=1}^{\infty}\frac{q^{k(k-1)/2}(-aq^{n})^{k}}{(q;q)_{k}}=-aq^{n}\sum_{k=0}^{\infty}\frac{q^{k(k+1)/2}(-aq^{n})^{k}}{(q;q)_{k+1}}.\]
Hence for $a>0$, we have \begin{eqnarray*}
\left|\sum_{k=0}^{\infty}\frac{q^{k(k+1)/2}(-aq^{n})^{k}}{(q;q)_{k+1}}\right| & \le & \frac{1}{1-q}\sum_{k=0}^{\infty}\frac{q^{k(k-1)/2}}{(q;q)_{k}}(aq^{n+1})^{k}\\
 & \le & \frac{(-aq^{2};q)_{\infty}}{1-q},\end{eqnarray*}
 and \eqref{eq:1.26} follows. Moreover from \eqref{eq:1.4} \begin{eqnarray*}
R_{2}(a;n) & = & aq^{n}\sum_{k=0}^{\infty}\frac{\left(aq^{n}\right)^{k}}{(q;q)_{k+1}},\end{eqnarray*}
 hence \begin{eqnarray*}
\left|R_{2}(a;n)\right| & \le & \frac{aq^{n}}{(1-q)(aq^{n};q)_{\infty}}\le\frac{aq^{n}}{(1-q)(aq;q)_{\infty}},\end{eqnarray*}
 which is \eqref{eq:1.28} and the proof of the lemma is complete.
\end{proof}

\section{Ramanujan's Entire Function $A_{q}(z)$}

For Ramanujan's entire function $A_{q}(z)$, we have the following:

\begin{thm}
\label{thm:q_airy}Given an arbitrary nonzero complex number $u$,
we have
\begin{enumerate}
\item For any positive rational number $t$ and $\lambda\in\mathbb{S}(t)$,
there are infinitely many positive integers $n$ and $m$ such that
\begin{equation}
tn=m+\lambda.\label{eq:2.1}\end{equation}
 For each such $\lambda$, $n$ and $m$ we have\begin{equation}
A_{q}(q^{-nt}u)=\frac{(-u)^{\left\lfloor m/2\right\rfloor }\left\{ \theta(-u^{-1}q^{\chi(m)+\lambda};q^{2})+r(n)\right\} }{(q;q)_{\infty}q^{\left\lfloor m/2\right\rfloor (nt-\left\lfloor m/2\right\rfloor )}}\label{eq:2.2}\end{equation}
 with\begin{gather}
|r(n)|\le\frac{3(-q^{3};q)_{\infty}\theta(\left|u\right|^{-1};q)}{1-q}\nonumber \\
\times\left\{ q^{m/4}+\frac{q^{m^{2}/16}}{|u|^{\left\lfloor m/4\right\rfloor +1}}\right\} .\label{eq:2.3}\end{gather}
 for $n$ sufficiently large.
\item For any positive irrational number $t$ and $\beta\in[0,1)$, there
are infinitely many positive integers $n$ and $m$ such that \begin{equation}
nt=m+\beta+\gamma_{n}\label{eq:2.4}\end{equation}
 with\begin{equation}
|\gamma_{n}|\le\frac{3}{n}.\label{eq:2.5}\end{equation}
 Let \begin{equation}
\nu_{n}=\left\lfloor -\frac{q^{2}\log n}{\log q}\right\rfloor \gg4,\label{eq:2.6}\end{equation}
for each such $\beta$, $n$ and $m$ we have\begin{equation}
A_{q}(q^{-nt}u)=\frac{(-u)^{\left\lfloor m/2\right\rfloor }\left\{ \theta(-u^{-1}q^{\chi(m)+\beta};q^{2})+e(n)\right\} }{(q;q)_{\infty}q^{\left\lfloor m/2\right\rfloor (nt-\left\lfloor m/2\right\rfloor )}}.\label{eq:2.7}\end{equation}
 and\begin{align}
|e(n)| & \le\frac{48(-q^{3};q)_{\infty}\theta(|u|^{-1};q)}{(1-q)}\nonumber \\
 & \times\left\{ \frac{\log n}{n}+q^{\nu_{n}^{2}/2}|u|^{\nu_{n}}+\frac{q^{\nu_{n}^{2}/2}}{|u|^{1+\nu_{n}}}\right\} .\label{eq:2.8}\end{align}
for $n$ sufficiently large.
\end{enumerate}
\end{thm}
\begin{proof}
From \eqref{eq:1.6}, we have \begin{equation}
A_{q}(q^{-nt}u)=\sum_{k=0}^{\infty}\frac{q^{k^{2}-knt}}{(q;q)_{k}}(-u)^{k}.\label{eq:2.9}\end{equation}
The classical Laplace method is used to study the asymptotics for\begin{equation}
\int_{-\infty}^{\infty}e^{\lambda f(x)}dx,\label{eq:2.10}\end{equation}
 as $\lambda\to+\infty$ , for example, see \cite{Wong}. If the real
function $f(x)$ has some maximas, from the nature of $e^{f(x)}$,
when $\lambda$ is large, then the major contributions of the integral
come from the neighbourhood of these maximas. We break the integral
into several pieces so that each piece has only one maxima and then
replace the integrands by simpler functions within each subintegrals
to get the asymptotics formula. Our situation is very similar here.
We notice that $q^{k^{2}-knt}$ has maximum around $\frac{nt}{2}$,
just as in the Laplace method for \eqref{eq:2.10}, we break the sum
into two subsums and estimate them respectively.

In the case that $t$ is a positive rational number, for any $\lambda\in\mathbb{S}(t)$
and $n$, $m$ are large, we have\begin{eqnarray}
A_{q}(q^{-nt}u)(q;q)_{\infty} & = & \sum_{k=0}^{\infty}(q^{k+1};q)_{\infty}q^{k^{2}-km-k\lambda}(-u)^{k}\nonumber \\
 & = & s_{1}+s_{2}\label{eq:2.11}\end{eqnarray}
 where\begin{equation}
s_{1}=\sum_{k=0}^{\left\lfloor m/2\right\rfloor }(q^{k+1};q)_{\infty}q^{k^{2}-km-k\lambda}(-u)^{k}\label{eq:2.12}\end{equation}
 and \begin{equation}
s_{2}=\sum_{k=\left\lfloor m/2\right\rfloor +1}^{\infty}(q^{k+1};q)_{\infty}q^{k^{2}-km-k\lambda}(-u)^{k}.\label{eq:2.13}\end{equation}
 In $s_{1}$ we reverse the order of summation to obtain\begin{eqnarray}
\frac{s_{1}q^{\left\lfloor m/2\right\rfloor (nt-\left\lfloor m/2\right\rfloor )}}{(-u)^{\left\lfloor m/2\right\rfloor }} & = & \sum_{k=0}^{\left\lfloor m/2\right\rfloor }(q^{\left\lfloor m/2\right\rfloor -k+1};q)_{\infty}q^{k^{2}}(-u^{-1}q^{\chi(m)+\lambda})^{k}\nonumber \\
 & = & \sum_{k=0}^{\infty}q^{k^{2}}(-u^{-1}q^{\chi(m)+\lambda})^{k}\nonumber \\
 & - & \sum_{k=\left\lfloor m/4\right\rfloor +1}^{\infty}q^{k^{2}}(-u^{-1}q^{\chi(m)+\lambda})^{k}\nonumber \\
 & + & \sum_{k=0}^{\left\lfloor m/4\right\rfloor }q^{k^{2}}(-u^{-1}q^{\chi(m)+\lambda})^{k}\left((q^{\left\lfloor m/2\right\rfloor -k+1};q)_{\infty}-1\right)\nonumber \\
 & + & \sum_{k=\left\lfloor m/4\right\rfloor +1}^{\left\lfloor m/2\right\rfloor }(q^{\left\lfloor m/2\right\rfloor -k+1};q)_{\infty}q^{k^{2}}(-u^{-1}q^{\chi(m)+\lambda})^{k}\nonumber \\
 & = & \sum_{k=0}^{\infty}q^{k^{2}}(-u^{-1}q^{\chi(m)+\lambda})^{k}+s_{11}+s_{12}+s_{13}.\label{eq:2.14}\end{eqnarray}
 Since \begin{equation}
0<(q^{\left\lfloor m/2\right\rfloor -k+1};q)_{\infty}<1\label{eq:2.15}\end{equation}
 for $0\le k\le\left\lfloor m/2\right\rfloor $, then,\begin{eqnarray}
|s_{11}+s_{13}| & \le & 2\sum_{k=\left\lfloor m/4\right\rfloor +1}^{\infty}q^{k^{2}}\left|u\right|^{-k}\nonumber \\
 & \le & 2\sum_{k=\left\lfloor m/4\right\rfloor +1}^{\infty}q^{k^{2}/2}\left|u\right|^{-k}\nonumber \\
 & \le & \frac{2q^{m^{2}/16}}{|u|^{\left\lfloor m/4\right\rfloor +1}}\theta(\left|u\right|^{-1};q).\label{eq:2.16}\end{eqnarray}
 By \eqref{eq:1.26}, for $0\le k\le\left\lfloor m/4\right\rfloor $,
we have\begin{equation}
\left|(q^{\left\lfloor m/2\right\rfloor -k+1};q)_{\infty}-1\right|\le\frac{(-q^{3};q)_{\infty}}{1-q}q^{m/4},\label{eq:2.17}\end{equation}
 then\begin{eqnarray}
|s_{12}| & \le & \frac{(-q^{3};q)_{\infty}q^{m/4}}{1-q}\sum_{k=0}^{\infty}q^{k^{2}}\left|u\right|^{-k}\nonumber \\
 & \le & \frac{(-q^{3};q)_{\infty}q^{m/4}}{1-q}\sum_{k=0}^{\infty}q^{k^{2}/2}\left|u\right|^{-k}\nonumber \\
 & \le & \frac{q^{m/4}(-q^{3};q)_{\infty}\theta(\left|u\right|^{-1};q)}{1-q},\label{eq:2.18}\end{eqnarray}
 hence\begin{equation}
\frac{s_{1}q^{\left\lfloor m/2\right\rfloor (nt-\left\lfloor m/2\right\rfloor )}}{(-u)^{\left\lfloor m/2\right\rfloor }}=\sum_{k=0}^{\infty}q^{k^{2}}(-u^{-1}q^{\chi(m)+\lambda})^{k}+r_{1}(n)\label{eq:2.19}\end{equation}
 with\begin{gather}
|r_{1}(n)|\le\frac{2(-q^{3};q)_{\infty}\theta(\left|u\right|^{-1};q)}{1-q}\nonumber \\
\times\left\{ q^{m/4}+\frac{q^{m^{2}/16}}{|u|^{\left\lfloor m/4\right\rfloor +1}}\right\} .\label{eq:2.20}\end{gather}
 In $s_{2}$ we change the summation from $k$ to $k+\left\lfloor m/2\right\rfloor $\begin{eqnarray}
\frac{s_{2}q^{\left\lfloor m/2\right\rfloor (nt-\left\lfloor m/2\right\rfloor )}}{(-u)^{\left\lfloor m/2\right\rfloor }} & = & \sum_{k=1}^{\infty}(q^{\left\lfloor m/2\right\rfloor +k+1};q)_{\infty}q^{k^{2}}(-uq^{-\chi(m)-\lambda})^{k}\nonumber \\
 & = & \sum_{k=1}^{\infty}q^{k^{2}}(-uq^{-\chi(m)-\lambda})^{k}\nonumber \\
 & + & \sum_{k=1}^{\infty}q^{k^{2}}(-uq^{-\chi(m)-\lambda})^{k}\left[(q^{\left\lfloor m/2\right\rfloor +k+1};q)_{\infty}-1\right]\nonumber \\
 & = & \sum_{k=-\infty}^{-1}q^{k^{2}}(-u^{-1}q^{\chi(m)+\lambda})^{k}+r_{2}(n).\label{eq:2.21}\end{eqnarray}
 By \eqref{eq:1.26}, for $k\ge1$\begin{equation}
\left|(q^{\left\lfloor m/2\right\rfloor +k+1};q)_{\infty}-1\right|\le\frac{q^{\left\lfloor m/2\right\rfloor +k+1}(-q^{3};q)_{\infty}}{1-q}\le\frac{(-q^{3};q)_{\infty}q^{m/2+k}}{1-q},\label{eq:2.22}\end{equation}
 then\begin{align}
|r_{2}(n)| & \le\frac{(-q^{3};q)_{\infty}q^{m/2}}{1-q}\sum_{k=1}^{\infty}q^{k^{2}}(|u|q^{1-\chi(m)-\lambda})^{k}\nonumber \\
 & \le\frac{(-q^{3};q)_{\infty}q^{m/2}}{1-q}\sum_{k=1}^{\infty}q^{k^{2}/2+k^{2}/2-k}|u|^{k}\nonumber \\
 & \le\frac{(-q^{3};q)_{\infty}q^{m/2-1/2}}{1-q}\sum_{k=1}^{\infty}q^{k^{2}/2}|u|^{k}\nonumber \\
 & =\frac{(-q^{3};q)_{\infty}q^{m/2-1/2}}{1-q}\sum_{k=-1}^{-\infty}q^{k^{2}/2}|u|^{-k}\nonumber \\
 & \le\frac{q^{m/4}(-q^{3};q)_{\infty}\theta(\left|u\right|^{-1};q)}{1-q}.\label{eq:2.23}\end{align}
 Thus we have proved that\begin{equation}
A_{q}(q^{-nt}u)=\frac{(-u)^{\left\lfloor m/2\right\rfloor }\left\{ \theta(-u^{-1}q^{\chi(m)+\lambda};q^{2})+r(n)\right\} }{(q;q)_{\infty}q^{\left\lfloor m/2\right\rfloor (nt-\left\lfloor m/2\right\rfloor )}}\label{eq:2.24}\end{equation}
 with\begin{eqnarray}
|r(n)| & \le & \frac{3(-q^{3};q)_{\infty}\theta(|u|^{-1};q)}{1-q}\nonumber \\
 & \times & \left\{ q^{m/4}+\frac{q^{m^{2}/16}}{|u|^{\left\lfloor m/4\right\rfloor +1}}\right\} .\label{eq:2.25}\end{eqnarray}
 for $n$, $m$ and $\lambda$ satisfying \eqref{eq:2.1} with $n$
and $m$ are sufficiently large.

In the case that $t$ is a positive irrational number, for any real
number $\beta\in[0,1)$, when $n$ and $m$ are sufficiently large
and satisfy \eqref{eq:2.4} and \eqref{eq:2.5} with \begin{equation}
1>\beta+\gamma_{n}>-1,\label{eq:2.26}\end{equation}
then\begin{equation}
2>\chi(m)+\beta+\gamma_{n}>-1.\label{eq:2.27}\end{equation}
For these integers $n$, we take\begin{equation}
\nu_{n}=\left\lfloor -\frac{q^{2}\log n}{\log q}\right\rfloor \gg4.\label{eq:2.28}\end{equation}
 Then, \begin{eqnarray}
A_{q}(q^{-nt}u)(q;q)_{\infty} & = & \sum_{k=0}^{\infty}(q^{k+1};q)_{\infty}q^{k^{2}-km-k\beta-k\gamma_{n}}(-u)^{k}\nonumber \\
 & = & s_{1}+s_{2},\label{eq:2.29}\end{eqnarray}
 with\begin{equation}
s_{1}=\sum_{k=0}^{\left\lfloor m/2\right\rfloor }(q^{k+1};q)_{\infty}q^{k^{2}-km-k\beta-k\gamma_{n}}(-u)^{k}\label{eq:2.30}\end{equation}
 and\begin{equation}
s_{2}=\sum_{k=\left\lfloor m/2\right\rfloor +1}^{\infty}(q^{k+1};q)_{\infty}q^{k^{2}-km-k\beta-k\gamma_{n}}(-u)^{k}.\label{eq:2.31}\end{equation}
 In $s_{1}$ we reverse the order of summation to get, \begin{eqnarray}
\frac{s_{1}q^{\left\lfloor m/2\right\rfloor (nt-\left\lfloor m/2\right\rfloor )}}{(-u)^{\left\lfloor m/2\right\rfloor }} & = & \sum_{k=0}^{\left\lfloor m/2\right\rfloor }(q^{\left\lfloor m/2\right\rfloor -k+1};q)_{\infty}q^{k^{2}}(-u^{-1}q^{\chi(m)+\beta+\gamma_{n}})^{k}\nonumber \\
 & = & \sum_{k=0}^{\infty}q^{k^{2}}(-u^{-1}q^{\chi(m)+\beta})^{k}\nonumber \\
 & - & \sum_{k=\nu_{n}+1}^{\infty}q^{k^{2}}(-u^{-1}q^{\chi(m)+\beta})^{k}\nonumber \\
 & + & \sum_{k=0}^{\nu_{n}}q^{k^{2}}(-u^{-1}q^{\chi(m)+\beta})^{k}\left(q^{k\gamma_{n}}-1\right)\nonumber \\
 & + & \sum_{k=0}^{\nu_{n}}q^{k^{2}}(-u^{-1}q^{\chi(m)+\beta+\gamma_{n}})^{k}\left\{ (q^{\left\lfloor m/2\right\rfloor -k+1};q)_{\infty}-1\right\} \nonumber \\
 & + & \sum_{k=\nu_{n}+1}^{\left\lfloor m/2\right\rfloor }q^{k^{2}}(-u^{-1}q^{\chi(m)+\beta+\gamma_{n}})^{k}(q^{\left\lfloor m/2\right\rfloor -k+1};q)_{\infty}\nonumber \\
 & = & \sum_{k=0}^{\infty}q^{k^{2}}(-u^{-1}q^{\chi(m)+\beta})^{k}+s_{11}+s_{12}+s_{13}+s_{14}.\label{eq:2.32}\end{eqnarray}
Since\begin{equation}
0<(q^{\left\lfloor m/2\right\rfloor -k+1};q)_{\infty}<1\label{eq:2.33}\end{equation}
 for $\nu_{n}+1\le k\le\left\lfloor m/2\right\rfloor $ and\begin{equation}
0<q^{\chi(m)+\beta+\gamma_{n}}\le q^{-1},\label{eq:2.34}\end{equation}
then,\begin{equation}
\left(\nu_{n}+1\right)^{2}-(\nu_{n}+1)>\frac{(\nu_{n}+1)^{2}}{2},\label{eq:2.35}\end{equation}
then,\begin{eqnarray}
|s_{11}+s_{14}| & \le & \sum_{k=\nu_{n}+1}^{\infty}q^{k^{2}}|u|^{-k}+\sum_{k=\nu_{n}+1}^{\infty}q^{k^{2}}|u|^{-k}q^{-k}\nonumber \\
 & \le & \frac{2q^{(\nu_{n}+1)^{2}-\nu_{n}-1}}{|u|^{\nu_{n}+1}}\sum_{k=0}^{\infty}q^{k^{2}}(|u|^{-1}q^{2\nu_{n}+1})^{k}\nonumber \\
 & \le & \frac{2q^{(\nu_{n}+1)^{2}-\nu_{n}-1}}{|u|^{\nu_{n}+1}}\sum_{k=0}^{\infty}q^{k^{2}/2}|u|^{-k}\nonumber \\
 & \le & \frac{2q^{\nu_{n}^{2}/2}}{|u|^{1+\nu_{n}}}\theta(|u|^{-1};q)\label{eq:2.36}\end{eqnarray}
 Since $0<q<1$ and $\lim_{n\to\infty}\log n/n=0$, there exists a
positive integer $N$ such that for $n\ge N$ \begin{equation}
n\ge\frac{3q^{2}\log n}{\log q^{-1}}.\label{eq:2.37}\end{equation}
Hence by \eqref{eq:1.23} and \eqref{eq:2.5}\begin{equation}
|q^{k\gamma_{n}}-1|\le\nu_{n}|\gamma_{n}|e^{\nu_{n}|\gamma_{n}|}\le\frac{3q^{2}e}{\log q^{-1}}\frac{\log n}{n}.\label{eq:2.38}\end{equation}
for $n\ge N$ and $0\le k\le\nu_{n}$. Thus \begin{align}
|s_{12}| & \le\frac{3eq^{2}}{\log q^{-1}}\frac{\log n}{n}\sum_{k=0}^{\infty}q^{k^{2}}|u|^{-k}\nonumber \\
 & \le\frac{3e}{\log q^{-1}}\frac{\log n}{n}\sum_{k=0}^{\infty}q^{k^{2}/2}|u|^{-k}\nonumber \\
 & \le\frac{12\theta(|u|^{-1};q)}{\log q^{-1}}\frac{\log n}{n}\label{eq:2.39}\end{align}
From \eqref{eq:1.26}, \begin{align}
|(q^{\left\lfloor m/2\right\rfloor -k+1};q)_{\infty}-1|q^{-k} & \le\frac{(-q^{3};q)_{\infty}q^{\left\lfloor m/2\right\rfloor -2k+1}}{1-q}\nonumber \\
 & \le\frac{(-q^{3};q)_{\infty}q^{m/4}}{1-q}\label{eq:2.40}\end{align}
 for $0\le k\le\nu_{n}$ and $n$ sufficiently large. Thus\begin{align}
|s_{13}| & \le\frac{(-q^{3};q)_{\infty}q^{m/4}}{1-q}\sum_{k=0}^{\nu_{n}}q^{k^{2}}|u|^{-k}\nonumber \\
 & \le\frac{(-q^{3};q)_{\infty}q^{m/4}}{1-q}\sum_{k=0}^{\nu_{n}}q^{k^{2}/2}|u|^{-k}\nonumber \\
 & \le\frac{(-q^{3};q)_{\infty}\theta(|u|^{-1};q)}{1-q}q^{m/4}\nonumber \\
 & \le\frac{(-q^{3};q)_{\infty}\theta(|u|^{-1};q)}{1-q}\frac{\log n}{n},\label{eq:2.41}\end{align}
for $n$ sufficiently large.\begin{equation}
\frac{s_{1}q^{\left\lfloor m/2\right\rfloor (nt-\left\lfloor m/2\right\rfloor )}}{(-u)^{\left\lfloor m/2\right\rfloor }}=\sum_{k=0}^{\infty}q^{k^{2}}(-u^{-1}q^{\chi(m)+\beta})^{k}+e_{1}(n)\label{eq:2.42}\end{equation}
 with\begin{align}
|e_{1}(n)| & \le\frac{24(-q^{3};q)_{\infty}\theta(|u|^{-1};q)}{1-q}\nonumber \\
 & \times\left\{ \frac{1-q+\log q^{-1}}{2\log q^{-1}}\frac{\log n}{n}+\frac{q^{\nu_{n}^{2}/2}}{|u|^{1+\nu_{n}}}\right\} \nonumber \\
 & \le\frac{24(-q^{3};q)_{\infty}\theta(|u|^{-1};q)}{1-q}\nonumber \\
 & \times\left\{ \frac{\log n}{n}+\frac{q^{\nu_{n}^{2}/2}}{|u|^{1+\nu_{n}}}\right\} .\label{eq:2.43}\end{align}
 In \eqref{eq:2.43}, we have used the inequality \eqref{eq:1.24}
for $x=1-q$ to show that $1-q<\log q^{-1}$. 

Similarly, in $s_{2}$ we change summation from $k$ to $k+\left\lfloor m/2\right\rfloor $,
\begin{align}
\frac{s_{2}q^{\left\lfloor m/2\right\rfloor (nt-\left\lfloor m/2\right\rfloor )}}{(-u)^{\left\lfloor m/2\right\rfloor }} & =\sum_{k=1}^{\infty}q^{k^{2}}(-uq^{-\chi(m)-\beta})^{k}\nonumber \\
 & -\sum_{k=\nu_{n}+1}^{\infty}q^{k^{2}}(-uq^{-\chi(m)-\beta})^{k}\nonumber \\
 & +\sum_{k=1}^{\nu_{n}}q^{k^{2}}(-uq^{-\chi(m)-\beta})^{k}(q^{-k\gamma_{n}}-1)\nonumber \\
 & +\sum_{k=1}^{\infty}q^{k^{2}}(-uq^{1-\chi(m)-\beta-\gamma_{n}})^{k}q^{-k}\left\{ (q^{\left\lfloor m/2\right\rfloor +k+1};q)_{\infty}-1\right\} \nonumber \\
 & +\sum_{k=\nu_{n}+1}^{\infty}q^{k^{2}}(-uq^{-\chi(m)-\beta-\gamma_{n}})^{k}\nonumber \\
 & =\sum_{k=-\infty}^{-1}q^{k^{2}}(-u^{-1}q^{\chi(m)+\beta})^{k}+s_{21}+s_{22}+s_{23}+s_{24}.\label{eq:2.44}\end{align}
 Just as we has done with the sum $s_{1}$, we can show that\begin{align}
|s_{21}+s_{24}| & \le2\sum_{k=\nu_{n}+1}^{\infty}q^{k^{2}-2k}|u|^{k}\nonumber \\
 & \le2q^{\nu_{n}^{2}/2}|u|^{\nu_{n}}\sum_{k=1}^{\infty}q^{k^{2}+2\nu_{n}k}|u|^{k}\nonumber \\
 & \le2q^{\nu_{n}^{2}/2}|u|^{\nu_{n}}\sum_{k=1}^{\infty}q^{k^{2}/2}|u|^{k}\nonumber \\
 & \le2q^{\nu_{n}^{2}/2}|u|^{\nu_{n}}\sum_{k=-1}^{-\infty}q^{k^{2}/2}|u|^{-k}\nonumber \\
 & \le2q^{\nu_{n}^{2}/2}|u|^{\nu_{n}}\theta(|u|^{-1};q),\label{eq:2.45}\end{align}
 \begin{align}
|s_{22}| & \le\sum_{k=1}^{\nu_{n}}q^{k^{2}/2+k^{2}/2-2k}|u|^{k}|q^{-k\gamma_{n}}-1|\nonumber \\
 & \le\frac{3q^{2}e}{\log q^{-1}}\frac{\log n}{n}q^{-2}\sum_{k=1}^{\nu_{n}}q^{k^{2}/2}|u|^{k}\nonumber \\
 & \le\frac{3e}{\log q^{-1}}\frac{\log n}{n}\sum_{k=-1}^{-\infty}q^{k^{2}/2}|u|^{-k}\nonumber \\
 & \le\frac{12\theta(|u|^{-1};q)}{\log q^{-1}}\frac{\log n}{n},\label{eq:2.46}\end{align}
\begin{align}
|s_{23}| & \le\frac{(-q^{3};q)_{\infty}q^{m/2}}{1-q}\sum_{k=1}^{\infty}q^{k^{2}/2+k^{2}/2-k}|u|^{k}\nonumber \\
 & \le\frac{(-q^{3};q)_{\infty}q^{m/2-1/2}}{1-q}\sum_{k=-1}^{-\infty}q^{k^{2}/2}|u|^{-k}\nonumber \\
 & \le\frac{(-q^{3};q)_{\infty}\theta(|u|^{-1};q)}{1-q}\frac{\log n}{n}\label{eq:2.47}\end{align}
for $n$ sufficiently large. Thus\begin{equation}
\frac{s_{2}q^{\left\lfloor m/2\right\rfloor (nt-\left\lfloor m/2\right\rfloor )}}{(-u)^{\left\lfloor m/2\right\rfloor }}=\sum_{k=-\infty}^{-1}q^{k^{2}}(-u^{-1}q^{\chi(m)+\beta})^{k}+e_{2}(n)\label{eq:2.48}\end{equation}
with\begin{align}
|e_{2}(n)| & \le\frac{24(-q^{3};q)_{\infty}\theta(|u|^{-1};q)}{(1-q)}\nonumber \\
 & \times\left\{ \frac{\log n}{n}+q^{\nu_{n}^{2}/2}|u|^{\nu_{n}}\right\} .\label{eq:2.49}\end{align}
 Thus we have\begin{equation}
A_{q}(q^{-nt}u)=\frac{(-u)^{\left\lfloor m/2\right\rfloor }\left\{ \theta(-u^{-1}q^{\chi(m)+\beta};q^{2})+e(n)\right\} }{(q;q)_{\infty}q^{\left\lfloor m/2\right\rfloor (nt-\left\lfloor m/2\right\rfloor )}}\label{eq:2.50}\end{equation}
 with\begin{equation}
e(n)=e_{1}(n)+e_{2}(n)\label{eq:2.51}\end{equation}
 and\begin{align}
|e(n)| & \le\frac{48(-q^{3};q)_{\infty}\theta(|u|^{-1};q)}{(1-q)}\nonumber \\
 & \times\left\{ \frac{\log n}{n}+q^{\nu_{n}^{2}/2}|u|^{\nu_{n}}+\frac{q^{\nu_{n}^{2}/2}}{|u|^{1+\nu_{n}}}\right\} .\label{eq:2.52}\end{align}
 for $n$ sufficiently large. 
\end{proof}
\begin{rem}
We have the following remarks on Theorem \ref{thm:q_airy}:
\begin{enumerate}
\item \eqref{eq:1.7}, \eqref{eq:2.2} and \eqref{eq:2.3} imply the trivial
formula for $\theta(z;q)$\begin{equation}
\theta(z;q)=zq^{1/2}\theta(zq;q).\label{eq:2.53}\end{equation}
 
\item We can rewrite \eqref{eq:2.2} and \eqref{eq:2.3} into the form\[
\theta(-u^{-1}q^{\chi(m)+\lambda};q^{2})=\frac{A_{q}(q^{-nt}u)(q;q)_{\infty}q^{\left\lfloor m/2\right\rfloor (nt-\left\lfloor m/2\right\rfloor )}}{(-u)^{\left\lfloor m/2\right\rfloor }}+o(1),\]
where the $o(1)$ is uniform for $q$ in any compact subset of $(0,1)$
and $u$ in any compact subset of $\mathbb{C}\backslash\left\{ 0\right\} $. 
\item We can rewrite the formulas \eqref{eq:2.7} and \eqref{eq:2.8} into
the form\[
\theta(-u^{-1}q^{\chi(m)+\beta};q^{2})=\frac{(q;q)_{\infty}A_{q}(q^{-nt}u)q^{\left\lfloor m/2\right\rfloor (nt-\left\lfloor m/2\right\rfloor )}}{(-u)^{\left\lfloor m/2\right\rfloor }}+o(1),\]
where the $o(1)$ is uniform for $q$ in any compact subset of $(0,1)$
and $u$ in any compact subset of $\mathbb{C}\backslash\left\{ 0\right\} $. 
\end{enumerate}
\end{rem}

\section{A Class of Entire Functions}

The phenomenon demonstrated with $A_{q}(z)$ in Theorem \ref{thm:q_airy}
is universal for a class of entire basic hypergeometric function of
type\begin{equation}
f(z)=\sum_{k=0}^{\infty}\frac{(a_{1},\dotsc,a_{r};q)_{k}q^{lk^{2}}}{(b_{1},\dotsc,b_{s};q)_{k}}z^{k},\label{eq:3.1}\end{equation}
 with $l>0$, where \begin{equation}
(a_{1},\dotsc,a_{r};q)_{k}=\prod_{j=1}^{r}(a_{j};q)_{k}.\label{eq:3.2}\end{equation}
 A confluent basic hypergeometric series is formally defined as\begin{equation}
_{m}\phi_{n}\left(\begin{array}{c|c}
\begin{array}{c}
a_{1},\dotsc,a_{m}\\
b_{1},\dots,b_{n}\end{array} & q,z\end{array}\right)=\sum_{k=0}^{\infty}\frac{(a_{1},\dotsc,a_{m};q)_{k}}{(q,b_{1},\dotsc,b_{n};q)_{k}}z^{k}\left(-q^{(k-1)/2}\right)^{k(m+1-n)},\label{eq:3.3}\end{equation}
 when $m+1-n>0$. It is clear that the function\begin{equation}
_{m}\phi_{n}\left(\begin{array}{c|c}
\begin{array}{c}
a_{1},\dotsc,a_{m}\\
b_{1},\dots,b_{n}\end{array} & q,-zq^{(m+1-n)/2}\end{array}\right)\label{eq:3.4}\end{equation}
 is of the form \eqref{eq:3.1} with \begin{equation}
l=\frac{m+1-n}{2}.\label{eq:3.5}\end{equation}
For a clean statement of the following theorem, we also define\begin{equation}
c(r,s;q):=\frac{(b_{1},\dotsc,b_{s};q)_{\infty}}{(a_{1},\dotsc,a_{r};q)_{\infty}}.\label{eq:3.6}\end{equation}

\begin{thm}
\label{thm:entire}Assume that \begin{equation}
0\le a_{1},\dotsc,a_{r},b_{1},\dotsc,b_{s}<1.\label{eq:3.7}\end{equation}
and \begin{equation}
u\in\mathbb{C}\backslash\left\{ 0\right\} .\label{eq:3.8}\end{equation}
We have
\begin{enumerate}
\item For any positive rational number $t$ and $\lambda\in\mathbb{S}(t)$,
there are infinitely many positive integers $n$ and $m$ such that
\begin{equation}
tn=m+\lambda.\label{eq:3.9}\end{equation}
 For each such $\lambda$, $n$ and $m$ we have\begin{equation}
f(q^{-lnt}u)=\frac{u^{\left\lfloor m/2\right\rfloor }\left\{ \theta(u^{-1}q^{l\chi(m)+l\lambda};q^{2l})+r(n)\right\} }{c(r,s;q)q^{l\left\lfloor m/2\right\rfloor (nt-\left\lfloor m/2\right\rfloor )}}\label{eq:3.10}\end{equation}
 with\begin{eqnarray}
|r(n)| & \le & \left(\frac{2}{1-q}\right)^{r+s+1}\frac{{\displaystyle \prod_{j=1}^{s}(-b_{j}q^{2};q)_{\infty}}\theta(\left|u\right|^{-1};q^{l})}{{\displaystyle \prod_{j=1}^{r}}(a_{j};q)_{\infty}}\nonumber \\
 & \times & \left\{ q^{m/4}+\frac{q^{lm^{2}/16}}{|u|^{\left\lfloor m/4\right\rfloor +1}}\right\} .\label{eq:3.11}\end{eqnarray}
 for $n$ sufficiently large.
\item For any positive irrational number $t$ and $\beta\in[0,1)$, there
are infinitely many positive integers $n$ and $m$ such that \begin{equation}
nt=m+\beta+\gamma_{n}\label{eq:3.12}\end{equation}
 with\begin{equation}
|\gamma_{n}|\le\frac{3}{n}.\label{eq:3.13}\end{equation}
 Let \begin{equation}
\nu_{n}=\left\lfloor -\frac{q^{2}\log n}{\log q}\right\rfloor \gg4,\label{eq:3.14}\end{equation}
for each such $\beta$, $n$ and $m$ we have\begin{equation}
f(q^{-lnt}u)=\frac{u^{\left\lfloor m/2\right\rfloor }\left\{ \theta(u^{-1}q^{l\chi(m)+l\beta};q^{2l})+e(n)\right\} }{c(r,s;q)q^{l\left\lfloor m/2\right\rfloor (nt-\left\lfloor m/2\right\rfloor )}},\label{eq:3.15}\end{equation}
 with\begin{align}
|e(n)| & \le\frac{48\prod_{j=1}^{s}(-b_{j}q^{2};q)_{\infty}\theta(|u|^{-1};q^{l})}{(1-q)\prod_{j=1}^{r}(a_{j};q)_{\infty}}\nonumber \\
 & \times\left\{ \frac{\log n}{n}+q^{l\nu_{n}^{2}/2}|u|^{\nu_{n}}+\frac{q^{{l\nu}_{n}^{2}/2}}{|u|^{1+\nu_{n}}}\right\} .\label{eq:3.16}\end{align}
 for $n$ sufficiently large. \\

\end{enumerate}
\end{thm}
\begin{proof}
The function\begin{align}
f(q^{-lnt}u)c(r,s;q) & =\sum_{k=0}^{\infty}\prod_{j=1}^{s}\left\{ {1+R}_{1}(b_{j};k)\right\} \prod_{j=1}^{r}\left\{ 1+R_{2}(a_{j};k)\right\} q^{l(k^{2}-ntk)}u^{k}\nonumber \\
 & =s_{1}+s_{2},\label{eq:3.17}\end{align}
 where\begin{equation}
s_{1}=\sum_{k=0}^{\left\lfloor m/2\right\rfloor }\prod_{j=1}^{s}\left\{ {1+R}_{1}(b_{j};k)\right\} \prod_{j=1}^{r}\left\{ 1+R_{2}(a_{j};k)\right\} q^{l(k^{2}-ntk)}u^{k}\label{eq:3.18}\end{equation}
and\begin{equation}
s_{2}=\sum_{k=\left\lfloor m/2\right\rfloor +1}^{\infty}\prod_{j=1}^{s}\left\{ {1+R}_{1}(b_{j};k)\right\} \prod_{j=1}^{r}\left\{ 1+R_{2}(a_{j};k)\right\} q^{l(k^{2}-ntk)}u^{k}.\label{eq:3.19}\end{equation}
 By Lemma \ref{lem:1}, we have\begin{equation}
\left|\prod_{j=1}^{s}\left\{ {1+R}_{1}(b_{j};k)\right\} \prod_{j=1}^{r}\left\{ 1+R_{2}(a_{j};k)\right\} -1\right|\le\left(\frac{2}{1-q}\right)^{r+s}\frac{\prod_{j=1}^{s}(-b_{j}q^{2};q)_{\infty}q^{k}}{\prod_{j=1}^{r}(a_{j};q)_{\infty}},\label{eq:3.20}\end{equation}
 and\begin{equation}
\left|\prod_{j=1}^{s}\left\{ {1+R}_{1}(b_{j};k)\right\} \prod_{j=1}^{r}\left\{ 1+R_{2}(a_{j};k)\right\} \right|\le\frac{1}{\prod_{j=1}^{r}(a_{j};q)_{\infty}}.\label{eq:3.21}\end{equation}
The rest of the proof is very similar to the proof for Theorem \ref{thm:q_airy}.
\end{proof}
\begin{acknowledgement*}
The author thanks the referee for many excellent suggestions for this
paper. The referee also kindly pointed out that there exists a connection
between classical Airy function and $q$-Airy function. When I was
visiting professor Mourad E. H. Ismail this summer, he showed me a
formula of such kind, but I quickly found out that the formula was
false by letting the parameter $q$ decrease to $0$. I believe that
professor Ismail informed the authors of the formula.
\end{acknowledgement*}


\begin{thebibliography}{10}
\bibitem{AKhiezer}N. I. Akhiezer, \emph{The Classical Moment Problem
and Some Related Questions in Analysis}, English translation, Oliver
and Boyed, Edinburgh, 1965.

\bibitem{Andrews1}G. E. Andrews, q-series: Their development and
application in analysis, number theory, combinatorics, physics, and
computer algebra, CBMS Regional Conference Series, number 66, American
Mathematical Society, Providence, R.I. 1986.

\bibitem{Andrews2}G. E. Andrews, Ramanujan's \char`\"{}Lost\char`\"{}
Note book VIII: The entire Rogers-Ramanujan function, Advances in
Math. 191 (2005), 393--407.

\bibitem{Andrews3}G. E. Andrews, Ramanujan's \char`\"{}Lost\char`\"{}
Note book IX: The entire Rogers-Ramanujan function, Advances in Math.
191 (2005), 408--422.

\bibitem{Andrews4}G. E. Andrews, R. A. Askey, and R. Roy, Special
Functions, Cambridge University Press, Cambridge, 1999.

\bibitem{Deift1}P. Deift, Orthogonal Polynomials and Random Matrices:
a Riemann-Hilbert Approach, American Mathematical Society, Providence,
2000.

\bibitem{Deift2}P. Deift, T. Kriecherbauer, K. T-R. McLaughlin, S.
Venakides, and X. Zhou, Strong asymptotics of orthogonal polynomials
with respect to exponential weights, Comm. Pure Appl. Math. 52 (1999),
1491--1552.

\bibitem{Gasper}G. Gasper and M. Rahman, Basic Hypergeometric Series,
second edition Cambridge University Press, Cambridge, 2004.

\bibitem{Hayman}W. K. Hayman, On the zeros of a q-Bessel function,
Contemporary Mathematics, volume 382, American Mathematical Society,
Providence, 2005, 205--216.

\bibitem{Hua}Hua Loo Keng, \textit{Introduction to Number Theory},
Springer-Verlag, Berlin Heidelberg New York, 1982.

\bibitem{Ismail1}M. E. H. Ismail, Asymptotics of q-orthogonal polynomials
and a q-Airy function, Internat. Math. Res. Notices 2005 No 18 (2005),
1063--1088.

\bibitem{Ismail2}M. E. H. Ismail, Classical and Quantum Orthogonal
Polynomials in one Variable, Cambridge University Press, Cambridge,
2005.

\bibitem{Ismail3}M. E. H. Ismail and X. Li, Bounds for extreme zeros
of orthogonal polynomials, Proc. Amer. Math. Soc. 115 (1992), 131--140.

\bibitem{Ismail4}M. E. H. Ismail and D. R. Masson, q-Hermite polynomials,
biorthogonal rational functions, Trans. Amer. Math. Soc. 346 (1994),
63--116.

\bibitem{Ismail5}M. E. H. Ismail and C. Zhang, Zeros of entire functions
and a problem of Ramanujan, Advances in Math., (2007), to appear.

\bibitem{Ismail6}M. E. H. Ismail and R. Zhang, Scaled asymptotics
for q-polynomials, Comptes Rendus, submitted.

\bibitem{Ismail7}M. E. H. Ismail and R. Zhang, Chaotic and Periodic
Asymptotics for q-Orthogonal Polynomials, joint with Mourad E.H. Ismail,
International Mathematics Research Notices, accepted. 

\bibitem{Kajiwara}K. Kajiwara, T. Masuda, M. Noumi, Y. Ohta, Y. Yamada,
Hypergeometric solutions to the q-Painlev\textbackslash{}'\{e\} equations,
Internat. Math. Res. Notices 47 (2004), 2497--2521.

\bibitem{Koekoek}R. Koekoek and R. Swarttouw, The Askey-scheme of
hypergeometric orthogonal polynomials and its q-analogues, Reports
of the Faculty of Technical Mathematics and Informatics no. 98-17,
Delft University of Technology, Delft, 1998.

\bibitem{Mehta}M. L. Mehta, Random Matrices, third edition, Elsevier,
Amsterdam, 2004.

\bibitem{Ramanujan}S. Ramanujan, The Lost Notebook and Other Unpublished
Papers (Introduction by G. E. Andrews), Narosa, New Delhi, 1988.

\bibitem{Saff}E. B. Saff and V. Totik, Logarithmic Potentials With
External Fields, Springer-Verlag, New York, 1997.

\bibitem{Szego}G. Szeg\textbackslash{}H\{o\}, Orthogonal Polynomials,
Fourth Edition, Amer. Math. Soc., Providence, 1975.

\bibitem{Wong}R. Wong, Asymptotic Approximations of Integrals, Academic
Press, Boston, 1989.

\bibitem{Whattaker}E. T. Whittaker and G. N. Watson, A Course of
Modern Analysis, fourth edition, Cambridge University Press, Cambridge,
1927. 
\end{thebibliography}
\end{document}